\documentclass[a4paper,12pt]{article}
\usepackage[T1]{fontenc}
\usepackage{graphics}
\usepackage{amsmath}
\usepackage{url}
\usepackage{theorem}
\usepackage{amsfonts}
\usepackage{amssymb}
\usepackage{makeidx}
\usepackage{bbm}
\usepackage{vmargin}
\usepackage[usenames]{color}
\usepackage[T1]{fontenc}
\usepackage{enumerate}
\usepackage{hyperref}
\usepackage{color}
\newenvironment{prooff}{}{\hfill$\square$}
\numberwithin{equation}{section}

\newcommand{\R}{\mathbb{R}^N}
\newcommand{\eps}{\varepsilon}
\newtheorem{theorem}{Theorem}[section]
\newtheorem{lemma}[theorem]{Lemma}

\newtheorem{rem}{Remark}[section]

 \usepackage{authblk}
\newcommand{\qed}{\nobreak \ifvmode \relax \else
      \ifdim\lastskip<1.5em \hskip-\lastskip
      \hskip1.5em plus0em minus0.5em \fi \nobreak
      \vrule height0.75em width0.5em depth0.25em\fi}
\title{Gradient estimate and  a Liouville theorem for  a $p$-Laplacian evolution equation with a gradient nonlinearity}
\author{Amal Attouchi\thanks{
E-mail: attouchi@math.univ-paris13.fr}}
\affil{
Universit\'e Paris 13, Sorbonne Paris Cit\'e, \\
Laboratoire Analyse, G\'eom\'etrie et Applications, CNRS (UMR 7539), \\
93430 Villetaneuse, France. 
}
\date{}
\begin{document}

\maketitle

\begin{abstract}
In this paper, we establish  a local gradient estimate for a $p$-Lpalacian equation with a fast growing  gradient nonlinearity. With
 this estimate, we can prove a parabolic Liouville theorem for ancient solutions satisfying some growth restriction near infinity. 
\end{abstract}

\noindent
{\it 2010 Mathematics Subject Classification}. {\scriptsize 35K65, 35K92, 35B45, 35B53
35K55}.\\
{\it  Key-words}.  p-Laplacian; local gradient estimate; parabolic Liouville-type
Theorem; gradient non-linearity
\vspace{1mm}
\hspace{.05in}

\section{Introduction and  main result}
In this paper, we are interested in qualitative properties of solutions of the non-linear degenerate parabolic
equation
\begin{eqnarray}\label{eqli}
 u_t-\Delta_p u=|\nabla u|^q, 
\end{eqnarray}
where 
$\Delta_p u= div(|\nabla u|^{p-2}\nabla u)$, $q>p-1>1$.\\
The kind of result we are going to prove are gradient estimates for local solutions in time-space, and a Liouville type theorem for ancient solutions.
In the last years, gradient estimates have played a key role in geometry and PDE since at least
the early work of Bernstein. Gradient a priori estimates are fundamental for elliptic and parabolic equations, leading to Harnack inequalities, Liouville
theorems, and compactness theorems for both linear and nonlinear PDE.
For the corresponding elliptic equation of  \eqref{eqli}, gradient estimates
were first considered by Lions \cite{Lio} for the linear diffusion case $p=2$. These estimates
were based upon the Bernstein technique. Recently for the  possibly degenerate elliptic equation with $q>p-1>0$, Bidaut-Véron, Huidobro, Véron \cite{hvv} obtained a priori universal gradient estimate for equations on a domain $\Omega$ of $\R$
and they extended their estimates to equations on
complete non compact manifolds satisfying a lower bound estimate on
the Ricci curvature. These estimates allowed them to derive some Liouville type theorems.

It  is  natural to look also for parabolic Liouville-type-theorems. In the linear diffusion case $p=2$ and for $q>1$,  Souplet and Zhang \cite{SZHAN} obtained local gradient estimate for locally upper bounded solution of \eqref{eqli}  ($u\leq M$) of the form
$$|\nabla u|\leq C(p, N, q) \left(t^{\frac{-1}{q}}+R^{-1}+ R^{\frac{-1}{q-1}}\right)(M+1-u)   \quad\text{in}\quad B(x_0, R)\times (0 ,T). $$
Relying on this estimate they proved that, under some growth  condition at infinity, ancient solutions  in the whole of $\mathbb{R}^N$ are constant. Motivated by their result, we generalize the gradient estimate and Liouville theorem to the case $1<p-1<q$. 
We also require that the solution is locally lower bounded. Using a Bernstein method, we have the following gradient estimate.
\begin{theorem}\label{gradest}
 Let $q > p-1>1$, $x_0 \in\mathbb{R}^N$ and  $R,T>0$. We set  $Q_{T,R}=B(x_0,R)\times(0,T)$.
Let $u$ be  a solution in  $L^{\infty}((0,T); W^{1, \infty}(B(x_0, R))$  of 
$$\partial_t u-\Delta_p u=|\nabla u|^q\qquad\rm{in}\quad Q_{T, R}.$$
 Suppose that  $|u|\leq M$ for some constant $M\geq 1$. Then,
\begin{equation}\label{gradmil}
 |\nabla u|\leq C(p, N, q) \left(t^{\frac{-1}{q}}+R^{-1}+ R^{\frac{-1}{q-p+1}}\right) M \quad\text{in}\quad Q_{T,\frac{R}{2}}.
\end{equation}

\end{theorem}
For the Cauchy-Dirichlet  problem associated to \eqref{eqli}, a gradient estimate involving the $W^{1, \infty}$  norm of the initial data has been obtained in \cite{AmalJDE,SZHAN}. In Theorem \ref{gradest} we only use the local $L^{\infty}$ norm of the solution but we get a weaker   estimate regarding the exponent on the distance to the boundary $R$.\\
Recently, for the singular diffusion case  $1<p< 2$ and  for $q=p$, F. Wang \cite{linfengeo}  established  gradient estimates similar to \eqref{gradmil} for  smooth, upper bounded, local solutions to  \eqref{eqli} on a closed manifolds or on  complete
noncompact Riemannian manifolds evolving under a Ricci flow. These estimates are of the form:
\begin{equation}\label{wangest}
\frac{|\nabla u|}{1-u}(x,t)\leq  C(N, p)\left(R^{-1}+t^{\frac{-1}{p}}+K^{\frac{2}{p}}+K\right)\quad\text{in}\quad Q_{T,\frac{R}{2}}.
\end{equation}
where $K>0$ is a constant related to the Ricci flow and the sectional curvature of the manifold. 
These estimates allowed to the author to provide some Harnack inequalities for positive solutions of the  following $p$-Laplace heat equation
\begin{equation}
|z|^{p-2} z_t=\Delta_p z.
\end{equation} 
The estimates \eqref{wangest} have been obtained by deriving an equation for $ w= |\nabla v|^p$,  $v=f^{-1}(-u)$ and $f(s)=e^{s/(p-1)}-1$. For $q>p>2$, we take a different auxiliary function $f$, adapted to the degenerate diffusion case and  to the fast growing gradient non-linearity. 

As an application of the  gradient estimate \eqref{gradmil}, we can state the following Liouville theorem for \eqref{eqli}.
\begin{theorem}\label{paraliou}
Assume that $q>p-1>1$ and let $\sigma=\min\left(1, \frac{1}{q-p+1}\right)$.
Assume that $u\in L^{\infty}_{loc}((- \infty, 0); W_{loc}^{1, \infty}(\R))$ is a weak solution of 
$$u_t-\Delta_p u=|\nabla u|^q, x\in\R, -\infty<t<0,$$
satisfying 
\begin{equation} \label{condgrow}
|u(x,t)|=o(|x|^{\sigma}+|t|^{\frac{1}{q}}), \quad \text{as}\quad |x|^{\sigma}+|t|^{\frac{1}{q}}\to\infty.
\end{equation}
Then $u$ is  constant.
\end{theorem}

\begin{rem}
The growth hypothesis \eqref{condgrow} is important (see the example of the function $u(x,t)=x_1+t$). However, we do not know if the exponents are sharp.
\end{rem}

Besides the works mentioned above, there are few  other studies on  gradient estimates and nonlinear Liouville theorems for a parabolic type equation on noncompact Riemannian manifolds. In this case the proof mostly relies on two types of gradient estimates or a combination of them. These estimates are known as Hamilton gradient estimate  (the estimate  only involves $\nabla u$ and $u$) \cite{hamil} and  Li-Yau's gradient estimate (the estimate involves $\nabla u, u$ and $u_t$) \cite{liyau}. Let us also mention that   the linear heat  equation on noncompact manifolds was studied by Souplet and Zhang in  \cite{souzhan} where they obtained a local gradient estimate related to the elliptic Cheng-Yau estimate and Hamilton's estimate  for the heat equation on compact manifolds. A Liouville theorem was also proved in \cite{souzhan}.  Hamilton-type gradient estimates were also used in \cite{meng, zhaoma, zhu}.
For $q=p>1$, a nonlinear analogue of Li-Yau's  estimate has been established in \cite{entroplap} for positive solutions of \eqref{eqli} on compact manifolds with nonnegative Ricci curvature. In \cite{entroplap}, the gradient estimate was not used to get Liouville theorems but to obtain an entropy formula. Nevertheless, Liouville theorems should be obtained as a consequence of the obtained gradient estimate.

This paper is organized as follows: In Section 2, we provide the proof of the  gradient estimate \eqref{gradmil}  and  we prove 
Theorem \ref{paraliou}.
In Sections 3 we give the proof of a technical auxiliary lemma that appears in the proof of the gradient estimate.

\section{Bernstein-type gradient estimate }
The proof of Theorem \ref{gradest} is based on the following technical lemma which is based on a Bernstein method. The most significant difficulty being the choice of the auxiliary function $f$ and the estimates coming from the cut-off argument. Let us mention that for  different suitable choice of $f$, gradient bounds global in space for the Cauchy problem  associated to \eqref{eqli} have been obtained in \cite{bartier}.\\
First let us make precise that by local weak solution of \eqref{eqli} we mean a function $u\in C_{loc}(\Omega\times (0,T))\cap L^{\infty}_{loc}(0,T; W_{loc}^{1,\infty}(\Omega))$  where  $\Omega$ is a smooth domain and 
such that   the integral equality
\begin{align*}
 \int_{\Omega}\left(u(x,t)\psi(x,t)-u(x,s)\psi(x,s)\right)\, dx &+\int_s^t\int_{\Omega} \left(-u \psi_t+|\nabla u|^{p-2} \nabla u\cdot\nabla \psi\right) \, dx\, d\tau\\
&=\int_s^t\int _{\Omega} |\nabla u|^q \psi \,dx \,d\tau
\end{align*}
holds for all $0 < s < t < T$ and for all for all  testing function $\psi\in C^1(\overline{\Omega}\times[0, T])$ such that $ \psi= 0$ near $\partial\Omega\times(0, T)$.\\
Now let  $\alpha \in  (0,1)$ to be chosen later on. Set $ R' =\frac{3R}{4}$. We select a cut-off function  $\eta$
$\in C^2(\bar{B}(x_o, R'))$, $0 \leq\eta \leq1$,  satisfying  $\eta = 0$ for $|x-x_0| = R'$  and such that 
\begin{eqnarray}
\left.
\begin{array}{rr}
|\nabla \eta|\leq C R^{-1} \eta ^{\alpha}\\
|D^2 \eta|+\eta^{-1}|\nabla \eta|^2\leq C R^{-2} \eta^{\alpha}
\end{array}
\right\}\,\,\,\mathrm{for} \, |x-x_0|<R',
\end{eqnarray}
for some $C=C(\alpha)>0$ (see \cite{souzhan} for  the existence of such  function).
\begin{lemma}\label{calcber}
 Assume that  $u$ is a local weak solution of \eqref{eqli} and  that $|u|\leq M$ in $Q_{T, R}$ for some $M>1$.
We consider a $C^3$ smooth increasing function $f$ satisfying $f''>0$, the following  differential equation 
\begin{equation}\label{difeq}
\left(\dfrac{f''}{f'}\right)'+(p-1)(1+N)\left(\dfrac{f''}{f'}\right)^2=0
\end{equation}
  and mapping $\left[0,3 \right]$ onto $[-M, M]$. Defining $v=f^{-1}(-u)$, we set  $w= |\nabla v|^2$ and $z=\eta w$. Then at any point where $|\nabla u|>0$,  $z$ satisfies the following differential inequality

\begin{align}\label{lem}
\mathcal{L}(z)&\leq- 2(q-1)(f')^{q-2} f'' w^{\frac{q+2}{2}}\eta + C(p,N) (f')^{p-2} R^{-2}\eta^{\alpha}w^{\frac{p}{2}}\nonumber\\
&+C(p, q) R^{-1}\eta^{\alpha}\left[w^{\frac{p+1}{2}}(f')^{p-3}f''+ w^{\frac{q+1}{2}}(f')^{q-1}\right]
\end{align}
where 
\begin{equation}
\mathcal{L}(z):=\partial_t z- \mathcal{A} z+ H \cdot\nabla z
\end{equation}
with $\mathcal{A}$ is given by \eqref{eqA} $\mathcal{H }$ is given by \eqref{eqH}.
\end{lemma}
The proof of lemma \ref{calcber} is postponed to the the next section.

\subsection*{Proof of Theorem \ref{gradest}}
\begin{prooff}
 Let  $u\in L_{loc}^{\infty}((0, \infty); W^{1, \infty}_{loc}(\Omega))$  be a local weak solution of \eqref{eqli}. Since $u$ and $\nabla u$ are locally bounded, using the result of Di Benedetto and Friedman \cite{fried,Dib}, we get that $\nabla u$
is a  locally $\mathrm{H\ddot{o}lder}$ continuous function. Thus  $z$ is a continuous function on $\overline{B(x_0, R')}\times [0, T]=\overline{Q}$, for  any $0<T$. Therefore, unless $z\equiv 0$ in $\overline{Q}$,  $z$  must reach a positive maximum at some point $(\hat{x},\hat{t})\in \overline{B(x_0, R')}\times [t_0, T]$.  
Since $z=0$ on $\partial B_{R'}\times [0, T]$, we deduce that $\hat{x}\in B_{R'}$.
Since $z(\hat{x}, \hat{t})>0$, we have that $|\nabla u|=f'(v)|\nabla v| >0$ and hence we can use Lemma \ref{calcber}.

Now let us  take $f(s)= M(s+1)^{\gamma}-2M$ where $\gamma$ is given by 
\begin{equation}\label{boncoef}
\gamma=\gamma(p,N)= \dfrac{(p-1)(N+1)+1}{(p-1)(N+1)}
\end{equation}
It is easy to see that $f$ satisfies the differential equation \eqref{difeq} and $f', f''>0$ and $f$ maps $\left[0, 3^{ \frac{1}{\gamma}}-1\right]$ onto $[-M, M]$. 
Let us also note that $\gamma\geq1$ and $\gamma-1\leq \dfrac{1}{p-1}\leq 1$.\\
By Lemma \ref{calcber} we get that, in a small neighbourhood $\tilde{Q}$ of $(\hat{x}, \hat{t})$, $z$ satisfies
\begin{align*}
\mathcal{L} z &\leq - 2(q-1)(f')^{q-2} f'' w^{\frac{q+2}{2}}\eta + C(p,N, \alpha) (f')^{p-2} R^{-2}\eta^{\alpha}w^{\frac{p}{2}}\\
&+C(p, q, \alpha) R^{-1}\eta^{\alpha}\left[w^{\frac{p+1}{2}}(f')^{p-3}f''+ w^{\frac{q+1}{2}}(f')^{q-1}\right].
\end{align*}
Hence 
\begin{align*}
(f')^{1-q}\mathcal{L} z &\leq - 2(q-1)\dfrac{ f''}{f'} w^{\frac{q+2}{2}}\eta + C(p,N, \alpha) (f')^{p-1-q} R^{-2}\eta^{\alpha}w^{\frac{p}{2}}\nonumber\\
&+C(p, q, \alpha) R^{-1}\eta^{\alpha}\left[w^{\frac{p+1}{2}}(f')^{p-q-1}\dfrac{f''}{f'}+ w^{\frac{q+1}{2}}\right].
\end{align*}
Since $v\in\left[0, (3)^{\frac{1}{\gamma}}-1\right]$, $\gamma, M\geq 1$, we have $1\leq v+1\leq (3)^{\frac{1}{\gamma}}\leq 3 $ and hence
\begin{equation}\label{gert}
\dfrac{1}{3(p-1)(N+1)}\leq\left(\dfrac{f''}{f'}\right)\leq \dfrac{1}{(p-1)(N+1)}\leq 1
\end{equation}
Using \eqref{gert} together with  the fact that $1\leq M\leq f'$ and  $p-q-1<0$,   we get that
\begin{align*}
(f')^{1-q}\mathcal{L}z&\leq -\dfrac{2(q-1)}{3(p-1)(N+1)} w^{\frac{q+2}{2}}\eta +C(N,p, \alpha) R^{-2} \eta^{\alpha}w^{\frac{p}{2}}\\
&+C(p, q, \alpha) R^{-1}\eta^{\alpha}\left[w^{\frac{p+1}{2}}
 + w^{\frac{q+1}{2}}\right].
\end{align*}
We take  $\alpha=\mathrm{max} \left(\frac{q+1}{q+2}, \frac{p+1}{q+2}\right)$.
Using the Young's inequality and recalling that $\eta\leq 1$,  then

\begin{itemize}
\item
for the conjugate exponents $r_1=\frac{q+2}{p}$, $s_1=\frac{q+2}{q-p+2}$ we have that
\begin{align*}
C(N,p, \alpha) R^{-2} \eta^{\alpha}w^{\frac{p}{2}}&=\eta^{\frac{p}{q+2}}w^{\frac{p}{2}}\,C(N,p, q,\alpha)\eta^{\alpha-p/(q+2)}R^{-2} \\   
&\leq \eps_1(N, p, q) \eta w^{\frac{q+2}{2}}+ C(N, p, q, \alpha) R^{\frac{-2(q+2)}{q-p+2}},
\end{align*}
\item
for the conjugate exponents $r_2=\frac{q+2}{p+1}$, $s_2=\frac{q+2}{q-p+1}$ we have that
\begin{align*}
C(N,p, q,\alpha) R^{-1} \eta^{\alpha}w^{\frac{p+1}{2}}&=\eta^{\frac{p+1}{q+2}}w^{\frac{p+1}{2}}\,C(N,p,q, \alpha) R^{-1} \eta^{\alpha-\frac{p+1}{q+2}} \\   
&\leq \eps_2 \eta w^{\frac{q+2}{2}}+ C(N, p, q,\alpha) R^{\frac{-(q+2)}{q-p+1}} 
\end{align*}
\item
and finally for the conjugate exponent $r_3=\frac{q+2}{q+1}$, $s_3=(q+2)$ we have that
\begin{align*}
C(N,p, q, \alpha) R^{-1} \eta^{\alpha}w^{\frac{q+1}{2}}&=\eta^{\frac{q+1}{q+2}}w^{\frac{q+1}{2}}\,C(N,p,q, \alpha) R^{-1} \eta^{\alpha-\frac{q+1}{q+2}} \\   
&\leq \eps_3 \eta w^{\frac{q+2}{2}}+ C(N, p, q, \alpha) R^{-(q+2)}.
\end{align*}
\end{itemize}
Choosing $\eps_i$ in such way that $\eps_1+\eps_2+\eps_3=\dfrac{1}{4} \dfrac{2(q-1)}{3(p-1)(N+1)}$, we get that
\begin{eqnarray}
 \begin{array}{ll}
(f')^{1-q}\mathcal{L}z&\leq -\dfrac{(q-1)}{2 (p-1)(N+1)}w^{\frac{q+2}{2}}\eta
+C(N, p, q, \alpha) R^{\frac{-2(q+2)}{q-p+2}}\\
&+C(N, p, q,\alpha) R^{\frac{-(q+2)}{q-p+1}} +
 C(N, p, q, \alpha) R^{-(q+2)}.
\end{array}
\end{eqnarray}
Using the fact that 
$$
\left\{
\begin{array}{ll}
\frac{1}{q-p+1}\leq \frac{2}{q-p+2}\leq 1 \,\,\, \mathrm{for}\, q\geq p,\\
1\leq\frac{2}{q-p+2}\leq\frac{1}{q-p+1} \,\,\, \mathrm{for}\, q\leq p,
\end{array}
\right.
$$
we have that
\begin{equation*}
(f')^{1-q}\mathcal{L}z\leq -\dfrac{(q-1)}{2 (p-1)(N+1)}w^{\frac{q+2}{2}}\eta +C(N, p, q, \alpha)\left[R^{-\frac{q+2}{q-p+1}}+  R^{-(q+2)}\right].
\end{equation*}
Setting
$$A=A(R,p,q, N):=C(N, p, q)\left( R^{-\frac{1}{q-p+1}}
+  R^{-1}\right)^2$$
and using that $(f')^{q-1}\geq M^{q-1}\geq 1$,
it follows that
\begin{equation}
\mathcal{L}z\leq- \dfrac{(q-1)}{4 (p-1)(N+1)}z^{\frac{q+2}{2}} \,\,\, \mathrm{in} \,\left\{ (x,t)\in Q_{T,R'};\,  z(x,t)\geq A \right\}.
\end{equation}
Next for $\lambda=\lambda(q, N, p)>0$ suitably chosen, the function  $\psi(t)=\lambda  t^{\frac{-2}{q}}$ satisfies $$\psi'(t)\geq -\dfrac{(q-1)}{4 (p-1)(N+1)} \psi^{\frac{q+2}{2}}.$$
Now for  $t_0\in(0,T)$ fixed, we define $\tilde{z}(t):=z(t+t_0, x)-\psi(t)$. It is easy to see that
$$\mathcal{L}\tilde{z}\leq 0 \,\,\, \mathrm{in}\, \left\{(x,t)\in Q_{T-t_0, R'};\,  \tilde{z}(x,t)\geq A\right\}.$$
Since  $\tilde{z}(t)\leq  0$ for  $t>0$ sufficiently small, we deduce from the maximum principle that 
$\tilde{z}(t)\leq  A$, i.e. $ z(x, t+t_0) \leq A+ \psi(t)$ in  $Q_{T-t_0, R'}$.

Finally using that  $z = \eta |\nabla v|^2$, letting  $t_0$ to 0, we get that
$$|\nabla v|\leq  C(N, p, q)(A + t^{\frac{-2}{q}})^{1/2}.$$
Using that 
$$v+1=\left(2-\frac{u}{M}\right)^{\frac{1}{\gamma}}\qquad\text{with} \left|\frac{u}{M}\right|\leq 1,$$
we get
$$\nabla v=\frac{-1}{\gamma M}\left(2-\frac{u}{M}\right)^{\frac{1-\gamma}{\gamma}}\nabla u.$$
It follows that
\begin{equation}
|\nabla u| \leq M \gamma|\nabla v|\leq C(N, p, q)(A + t^{\frac{-2}{q}})^{1/2} M \quad\textrm{in}\, Q_{T,\frac{R}{2}}.
\end{equation}
Here we used the fact that $\left(2-\frac{u}{M}\right)^{\frac{\gamma-1}{\gamma}}\leq 1$.\\
Hence we have
 $$|\nabla u| \leq C(N, p,q)\left(R^{-1}+R^{\frac{-1}{q-p+1}}+t^{\frac{-1}{q}}\right) M\quad\textrm{in}\, Q_{T,\frac{R}{2}}.$$
and the proof of Theorem \ref{gradest} is complete.
\end{prooff}

\subsection*{Proof of Theorem \ref{paraliou}}
\begin{prooff}
Fix $x_0\in\R$ and $t_0\in (-\infty, 0)$. Take $R\geq 1, T=R^{\sigma q}$ and set $Q=B(0, R)\times(0, T)$. Now we consider the function $U:=u(x+x_0, t+t_0-T)$. Using \eqref{condgrow}, we have that $|U|\leq M_{R}$ in $\overline{Q}$, where
$$M_{R}:= \underset{B(x_0, R)\times (t_0-T, t_0)}{\sup}\, |u|=o(T^{\frac{1}{q}}+R^{\sigma})=o(R^{\sigma}),\quad \textrm{as}\quad R\to\infty.$$
Applying Theorem \ref{gradest} to $U$ in $Q$, we get that
\begin{align*}
|\nabla u(x_0, t_0)|&=|\nabla U(0, T)|\leq C(N, p, q)R^{-\sigma}M_R
\end{align*}
and the conclusion follows by sending $R$ to $+\infty$.
\end{prooff}

\section{Proof of Lemma \ref{calcber} }
Our proof consists of three steps.
\subsection*{Step 1: computations}
Let $f$ be a $C^3$-function to be determined. We assume that $f', f''>0$. We put $v=f^{-1}(-u)$ and $w=|\nabla v|^2$.
By a  straightforward computation, we have that $v$ satisfies the following equation
\begin{align}
\partial_t v&=(f')^{p-2} w^{\frac{p-2}{2}}\left[\Delta v+(p-2) \frac{\left\langle D^2 v, \nabla v, \nabla v\right\rangle}{w}\right]+(p-1)(f')^{p-3}f'' w^{\frac{p}{2}}- (f')^{q-1} w^{\frac{q}{2}}\nonumber\\
=&(f')^{p-2} w^{\frac{p-2}{2}}\left[\Delta v+(p-2) \frac{\nabla w\cdot\nabla v}{2w}\right]+(p-1)(f')^{p-3}f'' w^{\frac{p}{2}}- (f')^{q-1} w^{\frac{q}{2}}.\label{eqv}
\end{align}
For $i=1,..., N$, we set $v_i=\frac{\partial v}{\partial x_i}$. In a neighbourhood $\tilde{Q}:=\omega\times (\tau_1, \tau_2)$ of any point $(\hat{x}, \hat{t})\in Q_{T,R}$ for which $|\nabla u|=f'(v)|\nabla v|>0$, the equation is uniformly parabolic and hence 
differentiating \eqref{eqv} with respect to $x_i$, we have
\begin{align}
\partial_t v_i=&(f')^{p-2} w^{\frac{p-2}{2}}\left[\Delta v_i+\dfrac{p-2}{2}\left( \dfrac{\nabla w_i\cdot\nabla v+\nabla w\cdot \nabla v_i}{w}-\dfrac{w_i \nabla w\cdot\nabla v}{w^2}\right)\right]\nonumber\\
&+(p-2)(f')^{p-3}f''v_i w^{\frac{p-2}{2}}\left[\Delta v+(p-2) \dfrac{\nabla w\cdot\nabla v}{2w}\right]\nonumber\\
&+\dfrac{p-2}{2} (f')^{p-2} w_i w^{\frac{p-4}{2}}\left[\Delta v+(p-2) \dfrac{\nabla w\cdot\nabla v}{2w}\right]\label{eqvi}\\
&+(p-1)((f')^{p-3}f'')' v_i w^{\frac{p}{2}}- (q-1)(f')^{q-2}f'' v_i w^{\frac{q}{2}}\nonumber\\
&+\dfrac{p(p-1)}{2 }(f')^{p-3}f''w_i w^{\frac{p-2}{2}}- \dfrac{q}{2}(f')^{q-1}w_i w^{\frac{q-2}{2}}.\nonumber
\end{align}
Here and in all the manuscript, the variable $v$ is omitted in the expression  of $f', f'', \left(\frac{f''}{f'}\right)'$, etc. The equalities are understood in a classical sense in  $\tilde{Q}$.
Multiplying \eqref{eqvi} by $2v_i$, summing over $i$ and using that
\begin{equation*}
\left\langle D^2 v, \nabla v, \nabla v\right\rangle=\dfrac{1}{2}\nabla w\cdot\nabla v,\qquad\qquad \Delta w =2\nabla v\cdot\nabla  \Delta v+2|D^2 v|^2,
\end{equation*}
\begin{equation*}
\sum_i 2 (\nabla v_i\cdot\nabla w) v_i=|\nabla w|^2,\qquad \sum_i (\nabla w_i\cdot\nabla v) v_i=\left\langle D^2w ,\nabla v, \nabla v\right\rangle,
\end{equation*}
we get that 
\begin{align}
\partial_t w&= |\nabla u|^{p-2}\Delta w+(p-2)|\nabla u|^{p-4} \left\langle D^2 w , \nabla u,\nabla u\right\rangle-2|\nabla u|^{p-2}|D^2 v|^2 \nonumber\\
&+(p-2) (f')^{p-2}w^{\frac{p-4}{2}}\Delta v \left(\nabla v\cdot\nabla w\right)
+\frac{(p-2)}{2}(f')^{p-2} w^{\frac{p-4}{2}}|\nabla w|^2\nonumber\\
&+\dfrac{(p-2)(p-4)}{2} (f')^{p-2} w^{\frac{p-6}{2}}\left(\nabla v\cdot\nabla w\right)^2 \label{eqw}\\
&- q(f')^{q-1}w^{\frac{q-2}{2}}\nabla w\cdot\nabla v +(p(p-1)+(p-2)^2)(f')^{p-3}f''w^{\frac{p-2}{2}}\nabla w\cdot\nabla v\nonumber\\
&+2\left[(p-1)((f')^{p-3}f'')' w^{\frac{p+2}{2}}- (q-1)(f')^{q-2}f''  w^{\frac{q+2}{2}}+(p-2)(f')^{p-3}f''w^{\frac{p}{2}}\Delta v\right].\nonumber
\end{align}
Here, when passing from \eqref{eqvi} to \eqref{eqw}, the terms have been transformed according to
$$\begin{array}{lllll}
L^1_{t1}\to \tilde{L}^1_{t1} + \tilde{L}^1_{t3},& L^1_{t2}\to \tilde{L}^1_{t2}, &L^1_{t3} \to \tilde{L}^2_{t2},&
\quad L^1_{t4} +L^3_{t2}\to \tilde{L}^3,\\
L^2_{t1}\to \tilde{L}^5_{t3}, &L^2_{t2}\to \tilde{L}^4_{t3},& &\\
L^3_{t1}\to \tilde{L}^2_{t1},& &&\\
L^4_{t1}\to\tilde{L}^5_{t1},&L^4_{t2}\to\tilde{L}^5_{t2},&&\\
L^5_{t1}\to\tilde{L}^4_{t2},&L^5_{t2}\to\tilde{L}^4_{t1},&&
\end{array}$$
(with obvious labeling).\\
Hence $w$ satisfies
\begin{equation*}
\partial_t w-\mathcal{A}(w)-\mathcal{H}\cdot \nabla w=-2|\nabla u|^{p-2} |D^2v|^2+\mathcal{N} (w)
\end{equation*}
where
\begin{align}
\mathcal{A}(w)&=|\nabla u|^{p-2}\Delta w+(p-2)|\nabla u|^{p-4} \left\langle D^2 w , \nabla u,\nabla u\right\rangle\label{eqA},\\
\mathcal{H}&=(p-2) (f')^{p-2}w^{\frac{p-4}{2}}\Delta v \nabla v+\frac{(p-2)}{2}(f')^{p-2}w^{\frac{p-4}{2}} \nabla w,
\nonumber\\
&+\dfrac{(p-2)(p-4)}{2} (f')^{p-2} w^{\frac{p-6}{2}}\left(\nabla v\cdot\nabla w\right)\nabla v - q(f')^{q-1}w^{\frac{q-2}{2}}\nabla v\label{eqH}\\
&+(p(p-1)+(p-2)^2)(f')^{p-3}f''w^{\frac{p-2}{2}}\nabla v  \nonumber\\
\mathcal{N}(w)&=2(p-1)((f')^{p-3}f'')' w^{\frac{p+2}{2}}- 2(q-1)(f')^{q-2}f''  w^{\frac{q+2}{2}}\nonumber\\
&+2(p-2)(f')^{p-2}\frac{f''}{f'}w^{\frac{p}{2}}\Delta v.
\end{align}
\subsection*{Step 2: equation for $z$ and useful estimates}
We set $z=\eta w$. Defining  the operator $$\mathcal{L}(z):=\partial_t z-\mathcal{A}(z)-\mathcal{H}\cdot\nabla z,$$
we have that 
\begin{align*}
\mathcal{L}z&=\eta \mathcal{L} w+ w\mathcal{L}\eta-2 |\nabla u|^{p-2}\nabla\eta\cdot\nabla w
-2(p-2) |\nabla u|^{p-4}(\nabla u\cdot\nabla \eta)( \nabla w\cdot\nabla u)\\
&=-2|\nabla u|^{p-2} \nabla\eta\cdot\nabla w
-2(p-2)|\nabla u|^{p-4}(\nabla u\cdot\nabla \eta)( \nabla w\cdot\nabla u)\\
&\quad +\eta\mathcal{N}w +w\mathcal{L}\eta -2|\nabla u|^{p-2}  |D^2v|^2\eta.
\end{align*}

\subsection*{Leading estimates}
Recalling that $f$ is increasing and that   $f''>0$, we get the following estimates.
\begin{enumerate}

\item Estimate of  $\eta\mathcal{N}w$\\
\begin{eqnarray}
\begin{array}{lll}
 \eta\mathcal{N}w&\leq 2(p-1)((f')^{p-3}f'')' w^{\frac{p+2}{2}}\eta- 2(q-1)(f')^{q-2} f'' w^{\frac{q+2}{2}}\eta\\
&+\frac{(f')^{p-2}}{2}w^{\frac{p-2}{2}}| D^2v|^2\eta+ 2N(p-1)^2 (f')^{p-2}\left(\dfrac{f''}{f'}\right)^2w^{\frac{p+2}{2}}\eta.
\end{array}
\end{eqnarray}
Here we used that  $$2(p-2)|\frac{f''}{f'} w\Delta v|\leq 2 N(p-1)^2 w^2\left(\dfrac{f''}{f'}\right)^2+\dfrac{|D^2 v|^2}{2}.$$
\item Estimate of $w\mathcal{L}(\eta)$
\begin{itemize}
\item Estimate of $w\mathcal{A}(\eta)$
\begin{equation}
|w\mathcal{A}(\eta)|\leq (f')^{p-2} w^{\frac{p}{2}}(\sqrt{N}+(p-2))|D^2 \eta|.
\end{equation}
\item Estimate of  $|w H\cdot\nabla \eta|$
\begin{eqnarray}
\begin{array}{ll}
 |w H\cdot\nabla\eta|&\leq
\underbrace{(f')^{p-2}w^{\frac{p-2}{2}}\left(C_1(p,N, \delta_1)\eta^{-1}|\nabla\eta|^2 w+\delta_1 [D^2v|^2\eta\right)}_{1} \\
&+
\underbrace{(f')^{p-2}w^{\frac{p-2}{2}}\left(C_2(p,N, \delta_2)\eta^{-1}|\nabla\eta|^2 w+\delta_2 [D^2v|^2\eta\right)}_{2}\\
&+ \underbrace{(f')^{p-2}w^{\frac{p-2}{2}}\left(C_3(p,N, \delta_3)\eta^{-1}|\nabla\eta|^2 w+\delta_3 [D^2v|^2\eta\right)}_{3} \\
&+ 2(p-1)^2(f')^{p-3}f''w^{\frac{p+1}{2}}|\nabla\eta|+ q (f')^{q-1}w^{\frac{q+1}{2}}|\nabla\eta|.
\end{array}
\end{eqnarray} 
\end{itemize}
(1) comes from an estimate via the Young's inequality of  $|(p-2)w \Delta v\nabla v\cdot\nabla\eta|$. Recalling that $\nabla w= \left(2 D^2v, \nabla v\right)$, (2) comes from an estimate of  
$\left|\frac{(p-2)}{2}w \nabla w\cdot\nabla\eta\right|$ and (3) come from an estimate of $\left|\frac{(p-2)(p-4)}{2}w\left(\nabla v\cdot\nabla w\right)(\nabla v\cdot\nabla\eta)\right|$.
\item  Estimate of $2|\nabla u|^{p-2}|\nabla\eta\cdot\nabla w|$.\\
Using the Young inequality, we have 
$$2|\nabla u|^{p-2}| \nabla\eta\cdot\nabla w|\leq(f')^{p-2}w^{\frac{p-2}{2}}\left[ C_4(p,N, \delta_4)\eta^{-1}|\nabla\eta|^2w+\delta_4 |D^2v|^2\eta \right].$$
\item Estimate of  $2(p-2)(\nabla u\cdot\nabla \eta) (\nabla w\cdot\nabla u)$
\begin{equation*}
|2(p-2)(\nabla u\cdot \nabla \eta)(\nabla w\cdot \nabla u)|\leq  (f')^{2}w\left[C_5(N,p, \delta_5)\eta^{-1}|\nabla\eta|^2 w+
|D^ 2v|^2\eta\right].
\end{equation*}

\end{enumerate}
 Finally recalling that $\nabla u=f'\nabla v$ and  choosing  $\delta_i$ in such way that  $-2+\delta_1+\delta_2+\delta_3+\delta_4+\delta_5=-1$ and then recalling the properties of the function $\eta$, we arrive at
\begin{align*}
\mathcal{L}(z)&\leq 2(p-1)\eta\left[\left((f')^{p-3}f''\right)' w^{\frac{p+2}{2}}+ N(p-1) (f')^{p-2}\left(\dfrac{f''}{f'}\right)^2w^{\frac{p+2}{2}}\right]\\
&- 2(q-1)(f')^{q-2} f'' w^{\frac{q+2}{2}}\eta + C(p,N, \alpha) (f')^{p-2} R^{-2}w^{\frac{p}{2}}\eta^{\alpha}\\
&+C(p, q, \alpha) \eta^{\alpha} R^{-1}\left[w^{\frac{p+1}{2}}(f')^{p-3}f''+ w^{\frac{q+1}{2}}(f')^{q-1}\right].
\end{align*}
\subsection*{Step 3: suitable choice for the function $f$}
To get rid of the term 
\begin{align}
&\left((f')^{p-3}f''\right)' w^{\frac{p+2}{2}}+ N(p-1) (f')^{p-2}\left(\dfrac{f''}{f'}\right)^2w^{\frac{p+2}{2}}\\
&=(f')^{p-2}w^{\frac{p+2}{2}}\left[\left(\dfrac{f''}{f'}\right)'+(p-2)\left(\dfrac{f''}{f'}\right)^2 +(p-1)N\left(\dfrac{f''}{f'}\right)^2\right]\\
&\leq (f')^{p-2}w^{\frac{p+2}{2}}\left[\left(\dfrac{f''}{f'}\right)'+(p-1)(N+1)\left(\dfrac{f''}{f'}\right)^2\right].
\end{align}

we shall take a function $f$ satisfying the following differential equation
\begin{equation}
\left(\dfrac{f''}{f'}\right)'+(p-1)(1+N)\left(\dfrac{f''}{f'}\right)^2=0.
\end{equation}

Hence we get that
\begin{align*}
\mathcal{L}(z)&\leq- 2(q-1)(f')^{q-2} f'' w^{\frac{q+2}{2}}\eta + C(p,N, \alpha) (f')^{p-2} R^{-2}\eta^{\alpha}w^{\frac{p}{2}}\\
&+C(p, q,\alpha) R^{-1}\eta^{\alpha}\left[w^{\frac{p+1}{2}}(f')^{p-3}f''+ w^{\frac{q+1}{2}}(f')^{q-1}\right].
\end{align*}


\begin{thebibliography}{10}

\bibitem{benar} L. Amour, M. Ben-Artzi, \emph{
Global existence and decay for a viscous Hamilton-Jacobi equation},
Nonlinear Analysis, \textbf{31} (1998), 621--628.

\bibitem{AmalJDE}
A.~Attouchi, \emph{{Well-posedness and gradient blow-up estimate near the boundary for a Hamilton-Jacobi equation with degenerate diffusion}},
J. Differential Equations \textbf{253} (2012), 2474--2492.
\bibitem{bartier}
{{\relax J.-{Ph}}} Bartier,  {\relax{Ph.}}~Lauren{\c{c}}ot, \emph{{Gradient
  estimates for a degenerate parabolic equation with gradient absorption and
  applications}}, Journal of Functional Analysis \textbf{254} (2008), 851--878.
\bibitem{weis}
M. Ben-Artzi, Ph. Souplet, F.B. Weissler,
\emph{The local theory for viscous Hamilton-Jacobi equations in Lebesgue spaces},
J. Math. Pure. Appl, \textbf{81} (2002),  343--378
\bibitem{hvv}
 M.F  Bidaut-Veron, M. Garcia-Huidobro, L. Veron,  \emph{
{Local and global properties of solutions of quasilinear Hamilton-Jacobi equations}}, preprint, (2014).


\bibitem{friedben}
E. DiBenedetto, A. Friedman,\emph{
Regularity of solutions of nonlinear degenerate parabolic systems}, J. Reine.
Angew. Math.
\textbf{349}
(1984), 83--128.

\bibitem{fried}E. DiBenedetto, A. Friedman,
\emph{H\"older estimates for nonlinear degenerate parabolic systems}, J. Reine
Angew. Math.
\textbf{357}
(1985), 1--22.
\bibitem{Dib}
E. DiBenedetto,\emph{Degenerate Parabolic Equations},
Springer Verlag, Series
Universitext, New York, (1993).

\bibitem{hamil}
R.S. Hamilton, \emph{A matrix Harnack estimate for the heat equation}, Comm. Anal. Geom., \textbf{1} (1993), 113-126.
\bibitem{entroplap}
B. Kotschwar, L. Ni, \emph{Local gradient estimates of $p$-harmonic functions, $1/H$-flow and an entropy formula}, Ann. Sci,Ec. Norm. Sup\'er.,(4)  \textbf{42} (2009), 1-36.
\bibitem{lady}
O.A. Ladyzenskaja, V.A. Solonnikov,  N.N. Ural'ceva, \emph{{Linear and
  quasi-linear equations of parabolic type}}, Amer. Math. Society, 1968.
	
\bibitem{liyau}
P. Li, S.T. Yau, \emph{On the parabolic kernel of the Schr\''odinger operator}, Acta Math, \textbf{153} (1986), 153-201.
	

\bibitem{Lio} P.L. Lions, \textit{Quelques remarques sur les probl\`emes elliptiques quasilin\'eaires du second ordre}, J. Anal.Math. {\bf  45} (1985), 234--254.
\bibitem{zhaoma}
	L. Ma, L. Zhao,  X. Song, \emph{Gradient estimate for the degenerate parabolic equation $u_t=\Delta(F(u))+H(u)$
 on manifolds},
Journal of Differential Equations,
\textbf{244} (2008), 1157--1177.
\bibitem{pola}
P. Polacik, P. Quittner \emph{
A Liouville-type theorem and the decay of radial solutions of a semilinear heat equation},
Nonlinear Analysis-theory Methods and Applications, \textbf{64} (2006), 1679--1689. 
\bibitem{superlinear}
P.~Quittner,  Ph. Souplet,
  \emph{{Superlinear parabolic problems: Blow-up, global existence and steady states}},
Birkhauser Advanced Texts, 2007.
\bibitem{soulip}
Ph. Souplet, \emph{
An optimal Liouville-type theorem for radial entire solutions of the porous medium equation with source}, Journal of Differential Equations, \textbf{246} (2009),  3980-4005.


\bibitem{souzhan}
 Ph. Souplet, Qi S. Zhang, \emph{Sharp gradient
estimate and Yau's Liou
ville theorem for the
heat equation on noncompact manifolds}, Bull. London Math. Soc. \textbf{38} (2006), 1045-1053.


 \bibitem{SZHAN}
Ph. Souplet,  Qi S. Zhang, \emph{Global solutions of inhomogeneous Hamilton-Jacobi equations}, Journal d'analyse mathématique, \textbf{99} (2006), 355--395. 


\bibitem{linfengeo} Lin-Feng Wang, \emph{
Gradient estimates for the p-Laplace heat equation under the Ricci flow},
Advances in Geometry \textbf{13} (2013), 349-368.

\bibitem{meng} M. Wang,
\emph{Liouville theorems for the ancient solution of heat flows},
Prcoceeding of the American mathematical society
\textbf{139} (2011),  3491--3496
\bibitem{wu}
 J. Wu, \emph{Gradient estimates for a nonlinear diffusion equation on complete manifolds}, J. Partial
Differ. Equ. \textbf{23}, (2010), 68-79. 
 \bibitem{yan}
Y. Yang, \emph{Gradient estimates for a nonlinear parabolic equation on Riemannian manifolds},
Proc. Amer. Math. Soc,  \textbf{136}, (2008), 4095-4102. 
\bibitem{zhu}
X.Zhu, \emph{Hmilton's gradient estimate and Liouville theorems for porous medium equations on noncompact Riemannian manifolds},  J. Math. Anal. Appl., \textbf{402} (2013), 201-206.

	\end{thebibliography}
\end{document}